\documentclass{article}

\usepackage{arxiv}

\usepackage[utf8]{inputenc} 
\usepackage[T1]{fontenc}    
\usepackage{hyperref}       
\usepackage{url}            
\usepackage{booktabs}       
\usepackage{amsfonts}       
\usepackage{nicefrac}       
\usepackage{microtype}      
\usepackage{mathptmx}
\usepackage{mathrsfs}
\usepackage{lipsum}
\usepackage{tabularx}
\usepackage{subcaption}
\usepackage{graphicx}

\title{Transformer-based Koopman Autoencoder for Linearizing Fisher's Equation}

\author{
 Kanav Singh Rana \\
  School of Mathematical and Statistical Sciences\\
  Indian Institute of Technology Mandi\\
  175005, India \\
  \texttt{kanavrana6@gmail.com} \\
   \And
 Nitu Kumari \\
  School of Mathematical and Statistical Sciences\\
  Indian Institute of Technology Mandi\\
  175005, India \\
  \texttt{nitu@iitmandi.ac.in} \\
}

\begin{document}
\maketitle
\begin{abstract}
A Transformer-based Koopman autoencoder is proposed for linearizing Fisher's reaction-diffusion equation. The primary focus of this study is on using deep learning techniques to find complex spatiotemporal patterns in the reaction-diffusion system. The emphasis is on not just solving the equation but also transforming the system's dynamics into a more comprehensible, linear form. Global coordinate transformations are achieved through the autoencoder, which learns to capture the underlying dynamics by training on a dataset with 60,000 initial conditions. Extensive testing on multiple datasets was used to assess the efficacy of the proposed model, demonstrating its ability to accurately predict the system's evolution as well as to generalize. We provide a thorough comparison study, comparing our suggested design to a few other comparable methods using experiments on various PDEs, such as the Kuramoto-Sivashinsky equation and the Burger's equation. Results show improved accuracy, highlighting the capabilities of the Transformer-based Koopman autoencoder. The proposed architecture in is significantly ahead of other architectures, in terms of solving different types of PDEs using a single architecture. Our method relies entirely on the data, without requiring any knowledge of the underlying equations. This makes it applicable to even the datasets where the governing equations are not known.
\end{abstract}


\section{Introduction}
\label{Section: 1}

Understanding complex dynamical systems is crucial for understanding many different types of natural occurrences, particularly those governed by PDEs. PDEs offer a theoretical foundation for modeling of spatiotemporal systems in many different fields, including engineering, physics, and biology. Of all these equations, Fisher's reaction-diffusion equation is commonly used to describe the spatiotemporal evolution of interacting species \cite{adomian1995fisher}. The nonlinear and complex nature of these equations makes them difficult to solve \cite{wang1988exact}. Moreover, there is no general mathematical architecture that can be used to solve such nonlinear PDEs.
\par
Transformations aimed at linearizing systems are commonly associated with Koopman operator theory \cite{koopman1931hamiltonian}. Koopman operator theory, with its modern interpretation in dynamical systems theory \cite{doi:10.1137/21M1401243}, plays an important role in our approach. While explicit construction of Koopman operators has limitations, Dynamic Mode Decomposition (DMD) \cite{rowley2009spectral, RanaKumari+2023}, offers numerical approximations to these operators. Recent advancements, such as Extended Dynamic Mode Decomposition (EDMD) \cite{williams2015data}, involve lifting system observables into higher-dimensional spaces through nonlinear transformations. This helps in the study of nonlinear phenomena, but the manual selection of transformation functions can limit effectiveness. 
\par 
Current initiatives have shown the effectiveness of neural networks in constructing Koopman embeddings for system dynamics \cite{lusch2018deep, li2017extended, yeung2019learning}. This increased interest in neural networks, emphasize their utility in approximating Koopman eigenfunctions and eigenvalues. This work aims to provide coordinate transformations that are interpretable, even for complex and nonlinear systems. Our primary focus is on attaining models that accurately capture the underlying low-rank dynamics while reducing overfitting and preserving interpretability.
\par
In recent years, deep learning methodologies have emerged as powerful tools in diverse scientific domains \cite{rajendra2020modeling}. However, relying solely on deep neural networks may not be sufficient, as their models often lack interpretability and parsimony \cite{raissi2018deep}. Using the capabilities of deep neural networks 
\cite{vaswani2017attention}, this research attempts to address the solution of Fisher's reaction-diffusion equation via linearising them using Transformer-based Koopman autoencoder. We introduce a systematic approach to utilize neural networks to discover coordinate transformations that aids in linearizing Fisher's reaction-diffusion equation. 
\par  
To ensure the success of our neural network architecture, we identify critical components, including appropriate network design, handling of the identity transformation, and informed selection of spatiotemporal data for training \cite{gin2021deep}. Our architecture emerges as a robust method, capable of addressing various PDEs without necessitating alterations to its fundamental design. This adaptability is evident in our research where we conducted comparative analysis using two different datasets and architectures. Notably, our architecture exhibited superior performance across various datasets, showcasing its effectiveness in solving diverse PDEs.
\par 
In Section \ref{Section: 2}, we explore the details of our methodology, involving the foundational principles of autoencoders, Fisher's equation, and Koopman operator theory within dynamical systems. To properly train our models, we also go into great detail about how we obtain dataset. The findings of our analysis are provided in Section \ref{Section: 3}, along with a thorough comparison of our suggested approach's performance with alternative approaches. We conclude in Section \ref{Section: 4} with a detailed discussion and critical conclusions drawn from our research.

\section{Methodology}\label{Section: 2}
We shall examine the methodological strategy used in this part to unveil the complex dynamics of Fisher's equation. We will describe our methods in detail, incorporating Koopman operator theory and autoencoders into the domain of dynamical systems. We also describe our data acquisition process, which is essential to effectively educate and train our algorithms.
\begin{figure*}[t!]
	\begin{subfigure}{\columnwidth}
		\centering
		\includegraphics[width=0.79\columnwidth]{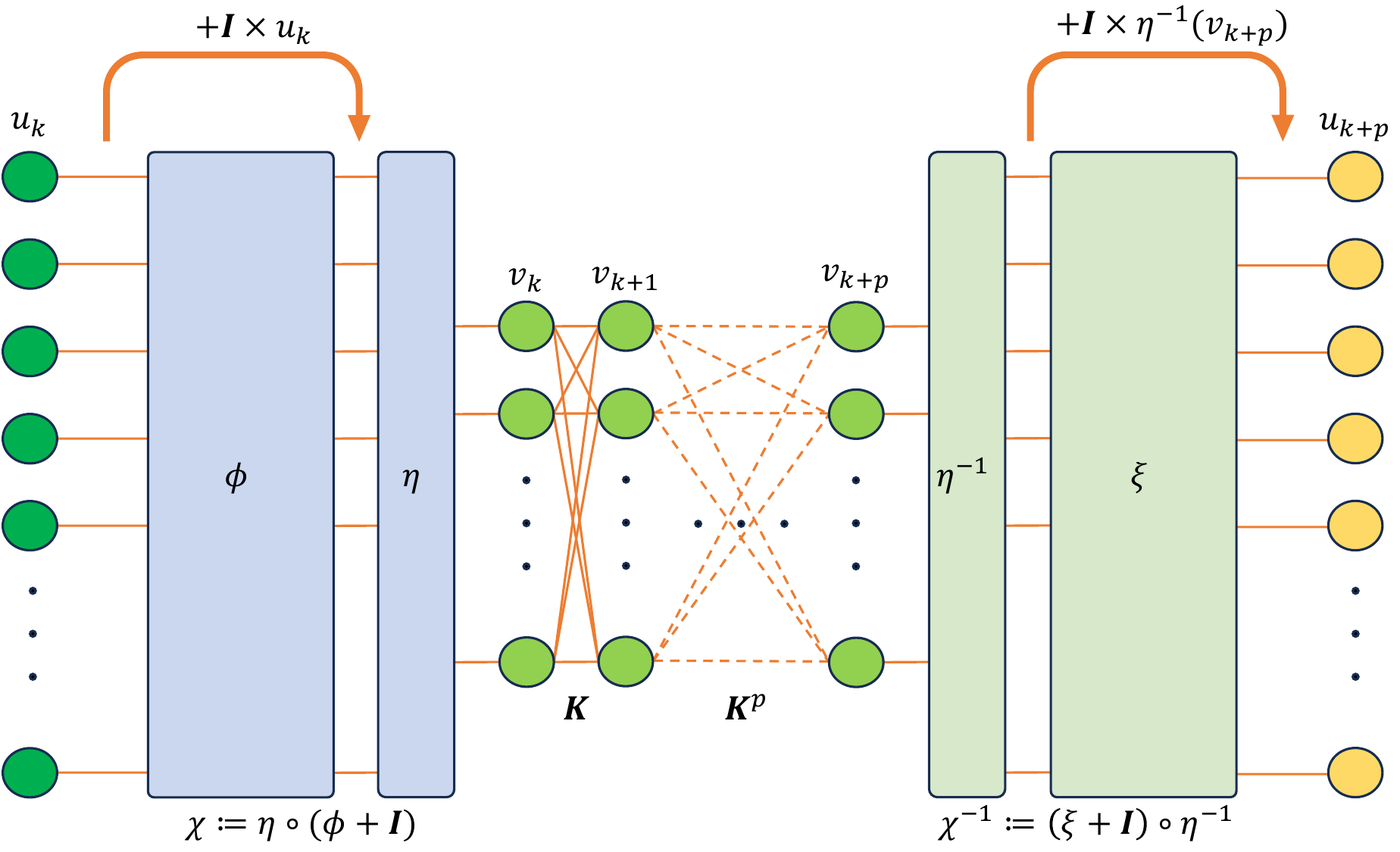}%
		\caption{}
		\label{Figure: 1_a}
	\end{subfigure}
	\begin{subfigure}{\columnwidth}
		\centering
		\includegraphics[width=0.80\columnwidth]{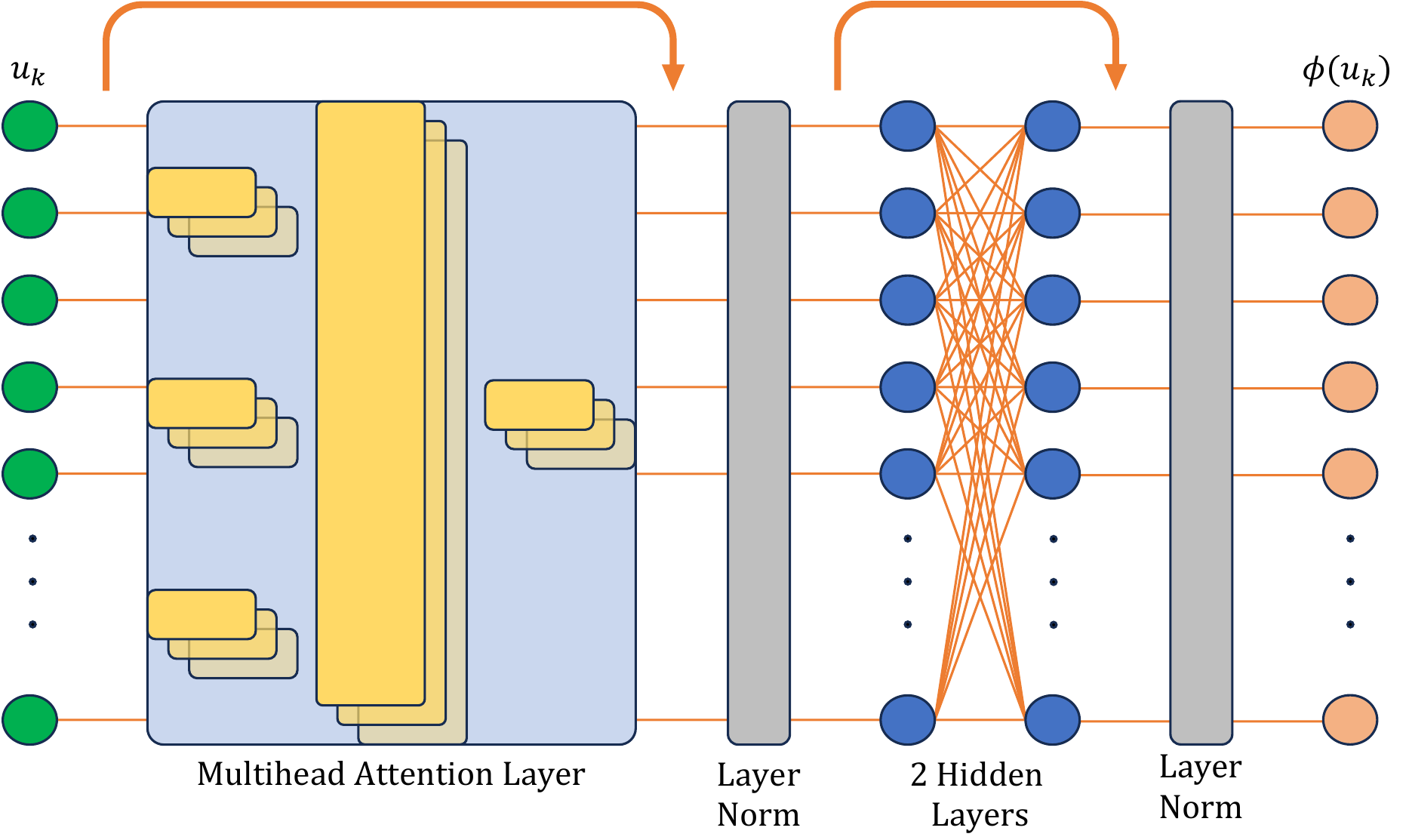}%
		\caption{}
		\label{Figure: 1_b}
	\end{subfigure}
	\caption{Neural network architecture. \textbf{(a)} Autoencoder comprise an outer encoder/decoder, an inner encoder/decoder and a matrix \textit{\textbf{K}}, which governs the dynamics of the system. \textbf{(b)} Transformer block that identifies the set of intrinsic coordinates of a linearized dynamical system.}
	\label{Figure: 1}
\end{figure*}
\subsection{Koopman Operator Theory} 
We examined nonlinear PDE of the form:
\begin{equation}
	u_{t}=\textit{N}(u,\;u_{x},\;u_{xx},\;t),
	\label{Eq: 1}
\end{equation}
where $u(x,t)\in \mathscr{M}\subseteq{\mathbb{R}}^n$ is the state of the system at time $t$ and space $x$ on a smooth manifold $\mathscr{M}$. Here $\textit{N}$ represents the governing equation that define the system.We will induce a mapping $\textbf{N}:\mathscr{M}\rightarrow \mathscr{M}$, where the state $u(x,t_{0})$ at time $t_{0}$ is mapped to a future time $t_{0} + t$. For partial differential equation, we may discretize the continuous function $u(x, t)$ at several spatial locations, $u(t)=[u(x_{1},t),\;u(x_{2},t),\; \dots,\;u(x_{n},t)]^{T} \in \mathbb{R}^n$. This yields discrete-time dynamical system \cite{kutz2013data}, for eq. (\ref{Eq: 1}):
\begin{equation}
	u_{k+1}=\textbf{N}(u_{k}), \quad 	k = 1,\;2,\dots,\;n.
\end{equation}
In the case where the function $\textbf{N}$ operates linearly, subsequent values of the state variable $u$ can be accurately predicted through a spectral expansion. However, due to the typical nonlinearity of $\textbf{N}$, a universally applicable scheme for solving such systems remains elusive.
\par 
Koopman operator theory presents an approach to linearizing nonlinear dynamical systems by means of the Koopman operator. This linear operator $\mathscr{K}$ is infinite-dimensional, analogous to $\textbf{N}$, such that, $\mathscr{K}g=g\circ \textbf{N}$, that advances the observable functions $g: \mathscr{M}\rightarrow \mathbb{C}$:
\begin{equation}
	\mathscr{K}g(u_{k}) = g(\textbf{N}(u_{k})) = g(u_{k+1}).
\end{equation}
The eigen-decomposition of operator $\mathscr{K}$ enables the simplification of the dynamics and representing solutions to the dynamical system. We seek to identify eigenfunctions $\varphi_k$ of the Koopman operator, corresponding to eigenvalues $\lambda_k$, satisfying:
\begin{equation}
	\varphi(u_{k+1})=\mathscr{K}\varphi(u_{k})=\lambda\varphi(u_k).
\end{equation}
\par 
Despite being an infinite-dimensional operator, the Koopman operator can be approximated in a finite-dimensional manner by considering the space defined by a finite set of Koopman eigenfunctions. When applied within this space, the Koopman operator can be converted into a matrix representation. In practice, procuring Koopman eigenfunctions may pose a greater challenge compared to obtaining the solution of eq. (\ref{Eq: 1}).

\subsection{Autoencoders}
Autoencoders are a type of artificial neural network used for unsupervised learning and dimensionality reduction tasks \cite{kramer1991nonlinear}. They consist of an encoder and a decoder, with the goal of learning a compressed representation of the input data. The encoder extracts the most crucial features from the input data and maps them to a lower-dimensional representation. This compressed representation, often referred to as the ``code" or ``latent variables," holds a condensed version of the input data.
\par
The entire architecture is depicted in figure \ref{Figure: 1_a}. The encoder in our design is $\chi$; the second component is the linear dynamics $\textit{\textbf{K}}$; the third component is the decoder $\chi^{-1}$. There are two segments in each of the encoder and decoder. The encoder is made up of the inner encoder $\eta$ and the outside encoder $\phi + \textit{\textbf{I}}$, where \textit{\textbf{I}} is an identity matrix. The inner encoder $\eta$ lowers dimensionality and diagonalizes the system, while the outside encoder conducts a coordinate transformation into a linear dynamics space. The inverses of the inner decoder $\eta^{-1}$ and outer decoder $\xi + \textit{\textbf{I}}$ are the inner and outer encoder, respectively. Both the outer encoder and outer decoder utilize a residual connection \cite{he2016deep}.
\par
The architecture function $\phi$, as shown in figure \ref{Figure: 1_b}, consists of a transformer block, comprising a multi-head self-attention mechanism followed by layer normalization and feedforward neural networks. This transformer block is applied to the input data $u_k$, allowing the model to capture intricate dependencies and patterns within the data. In addition, residual connections are employed within the block, enabling effective information flow through the network.
\par
The loss function utilized for training the network comprises five distinct components, each serving to enforce specific conditions. Loss function is given by: 
\begin{itemize}
	\item[\textbf{Loss 1:}] $\|u_{k}-\chi^{-1}(\chi(u_{k}))\|$ focuses on the autoencoder principle, aiming for a reversible transformation between the state space and intrinsic coordinates, where linear dynamics prevail. 
	\item[\textbf{Loss 2:}] $\|u_{k+p}-\chi^{-1}(\textit{\textbf{K}}^{p}\chi(u_{k}))\|$, ensures accurate forecasting of the next state given the current one, extending to multiple future time steps through iterative matrix multiplication. 
	\item[\textbf{Loss 3:}] $\|\chi(u_{k+p})-\textit{\textbf{K}}^{p}\chi(u_{k})\|$ emphasizes the requirement for linear dynamics within the intrinsic coordinates. \item[\textbf{Loss 4:}] $\|u_{k}-(\xi+\textit{\textbf{I}})((\phi+\textit{\textbf{I}})(u_{k}))\|$ addresses the outer autoencoder objective, emphasizing the separation of transformation and dimensionality reduction or diagonalization processes, allowing better interpretability. 
	\item[\textbf{Loss 5:}] $\|(\phi+\textit{\textbf{I}})u_{k}-\eta^{-1}(\eta((\phi+\textit{\textbf{I}})(u_{k})))\|$ mirrors \textbf{Loss 4} but pertains specifically to the outer encoder and decoder components, reinforcing the autoencoder principle.
\end{itemize} The total loss is incorporated by summing all the aforementioned losses, each carrying the same weight.

\subsection{Fisher's Equation}
For a state variable $u(x,t)$, we consider one-dimensional Fisher's equation:
\begin{equation}\label{eq: 5}
	\frac{\partial u(x,t)}{\partial t}=\alpha \frac{\partial^2 u(x,t)}{\partial x^2}+\beta u(x,t)(1-u(x,t)) \quad x\in\left(-\pi, \pi\right), t>0,
\end{equation}
with periodic boundary conditions, where $\alpha$ is the diffusion coefficient and $\beta>0$ is the reaction coefficient. For our results, we use $\alpha=1$ and $\beta=1$. We discretize the spatial domain, represented by $x$, into $n=128$ equidistant points.
\par 
The neural networks is trained using data generated through numerical simulations of the abovementioned PDE. We have curated a dataset comprising 60,000 initial conditions for training purposes. This PDE has been simulated forward for 50 time instances, with a time step $\Delta t=0.002$, culminating in a total duration of 0.10 seconds.
\par
Within the training dataset, we have incorporated three distinct types of initial conditions: white noise, sine waves, and square waves. Similarly, for validation and testing purposes, we have maintained consistency in terms of time and time steps. The validation set comprises 20,000 initial conditions, mirroring the types of initial conditions found in the training set. For the testing phase, we have assembled a dataset containing 21,000 initial conditions, encompassing seven different types. Three types of initial conditions in testing data are already present in the training data: white noise, sine waves, and square waves. Additionally, we have introduced Gaussian waves, triangle waves, sawtooth waves, and pulse waves, each comprising 3,000 initial conditions.
\begin{figure*}
	\begin{subfigure}{\columnwidth}
		\centering
		\includegraphics[width=0.70\columnwidth]{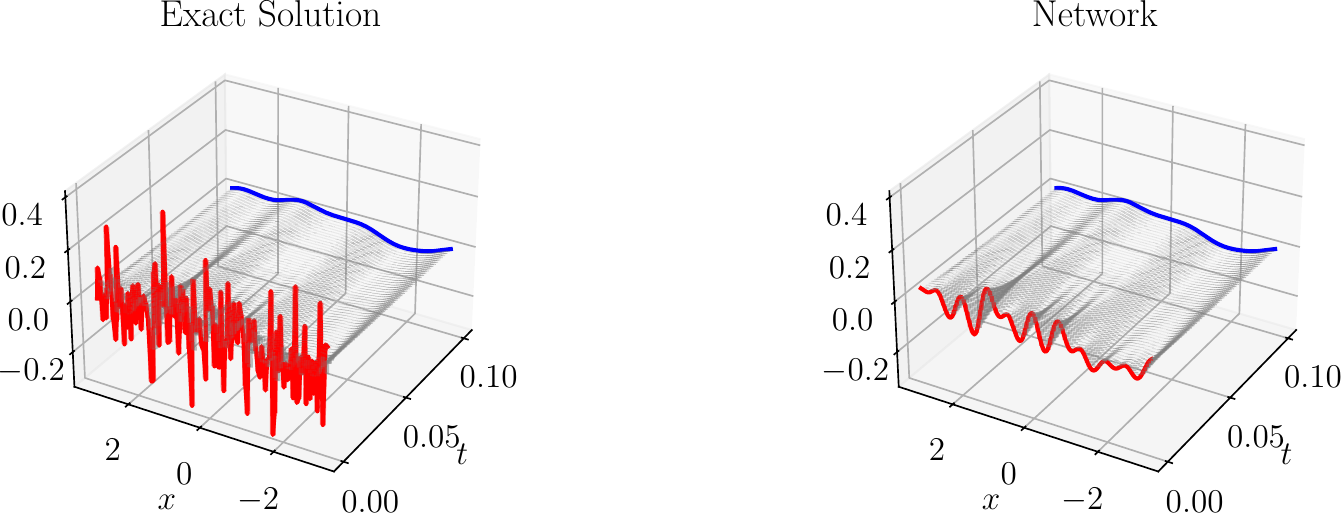}%
		\caption{}
		\label{Figure: 2_a}
	\end{subfigure}
	\begin{subfigure}{\columnwidth}
		\centering
		\includegraphics[width=0.70\columnwidth]{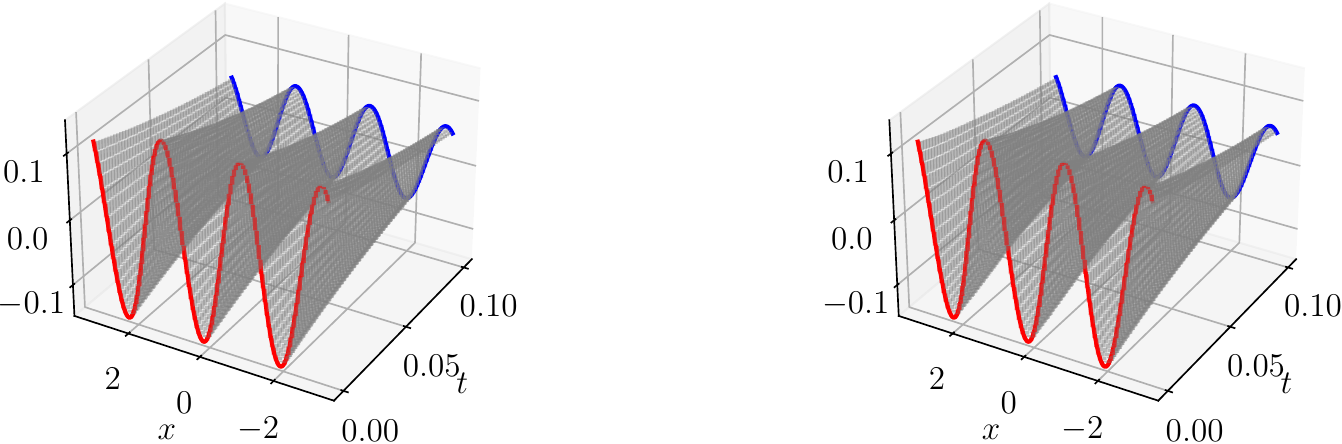}%
		\caption{}
		\label{Figure: 2_b}
	\end{subfigure}
	\begin{subfigure}{\columnwidth}
		\centering
		\includegraphics[width=0.70\columnwidth]{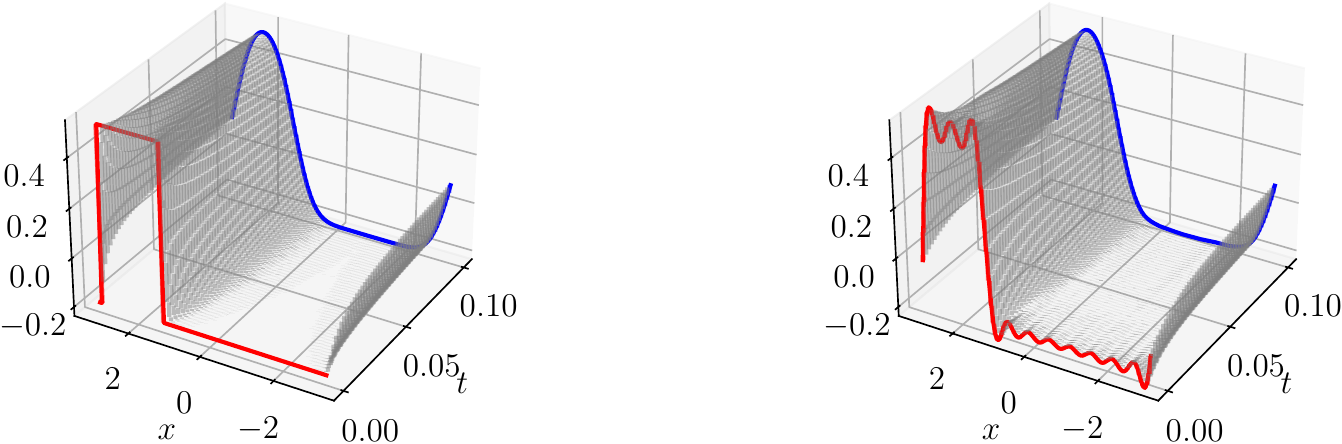}%
		\caption{}
		\label{Figure: 2_c}
	\end{subfigure}
	\caption{A comparison of exact solution and neural network predictions for Fisher's equation. The initial conditions utilized include \textbf{(a)} white noise, \textbf{(b)} sine waves, and \textbf{(c)} square waves, which are the same as the training dataset.}
	\label{Figure: 2}
\end{figure*}

\section{Results}\label{Section: 3}
\begin{figure*}
	\begin{subfigure}{\columnwidth}
		\centering
		\includegraphics[width=0.70\columnwidth]{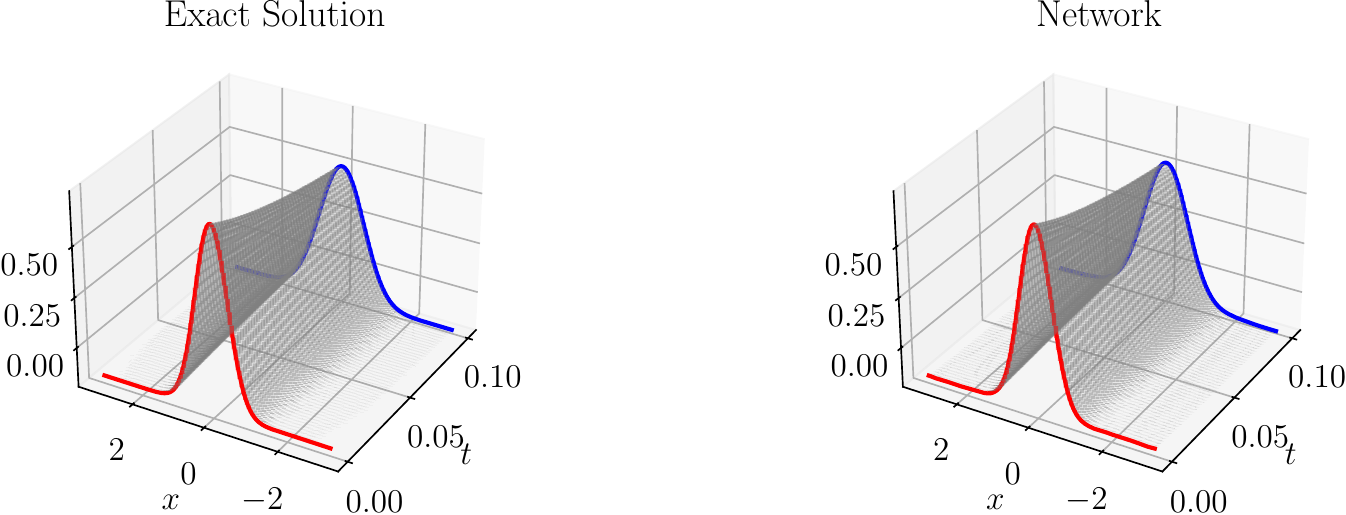}%
		\caption{}
		\label{Figure: 3_a}
	\end{subfigure}
	\begin{subfigure}{\columnwidth}
		\centering
		\includegraphics[width=0.70\columnwidth]{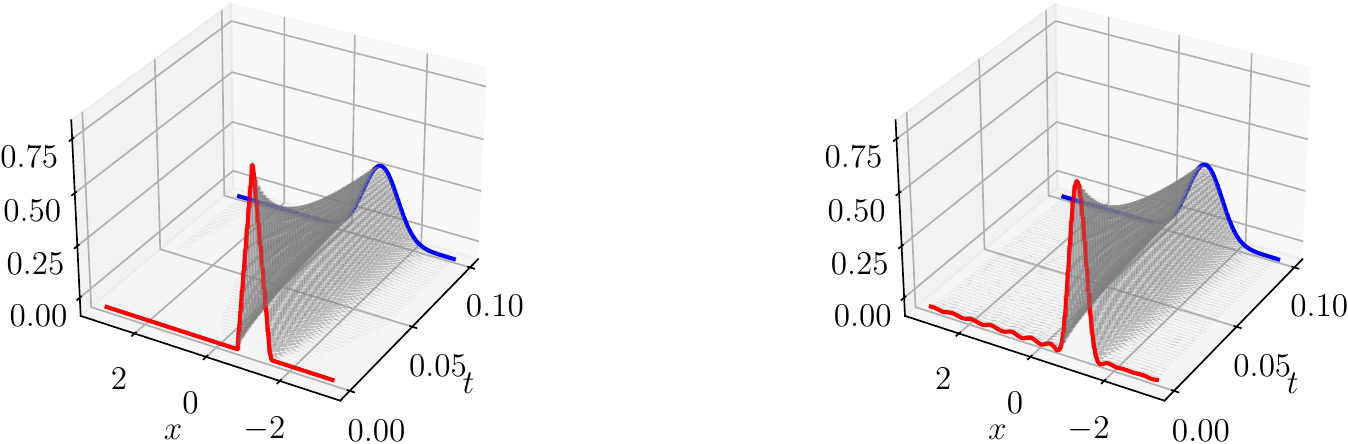}%
		\caption{}
		\label{Figure: 3_b}
	\end{subfigure}
	\begin{subfigure}{\columnwidth}
		\centering
		\includegraphics[width=0.70\columnwidth]{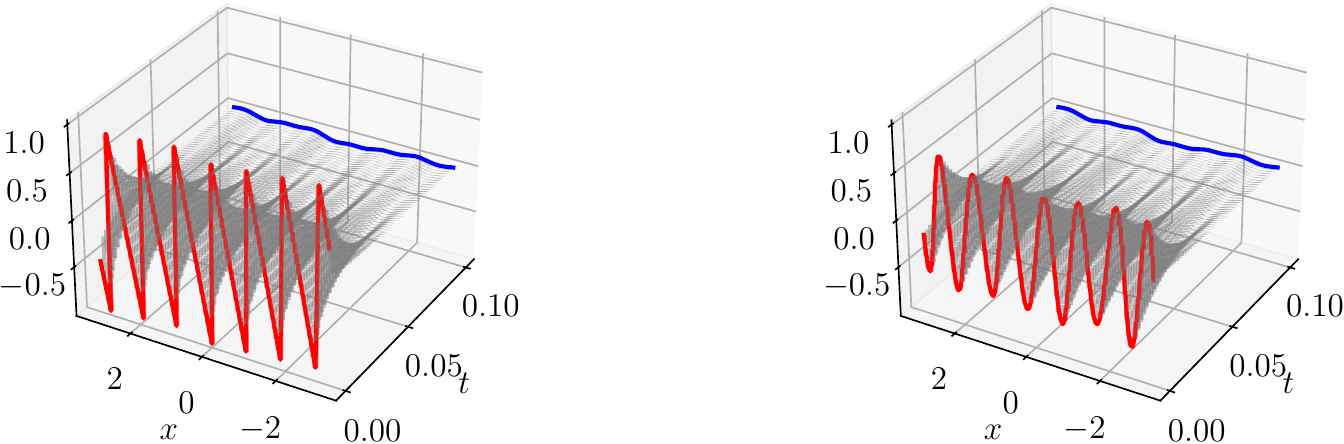}%
		\caption{}
		\label{Figure: 3_c}
	\end{subfigure}
	\begin{subfigure}{\columnwidth}
		\centering
		\includegraphics[width=0.70\columnwidth]{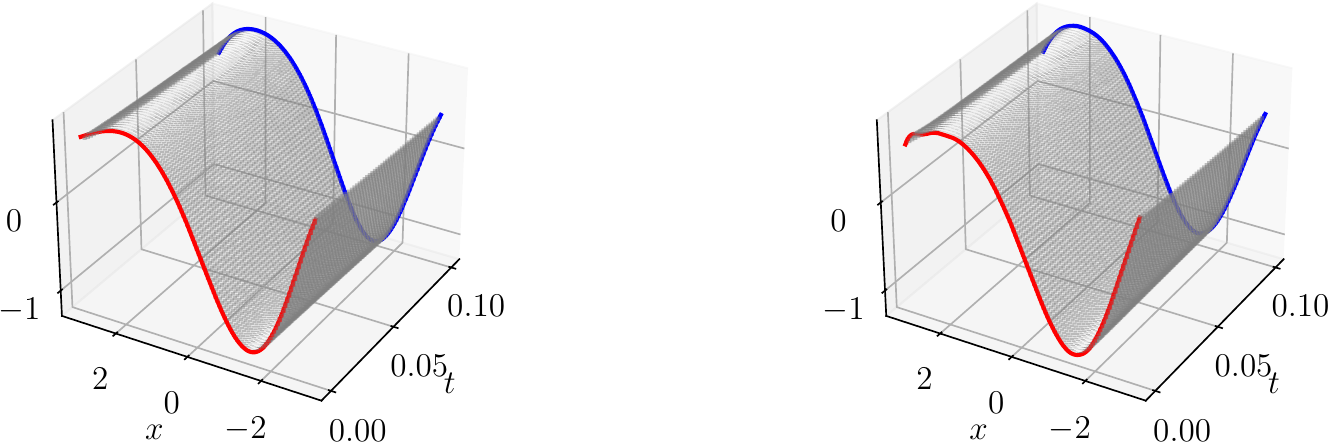}%
		\caption{}
		\label{Figure: 3_d}
	\end{subfigure}
	\caption{A comparison of exact solution and neural network predictions of Fisher's equation. The initial conditions utilized include \textbf{(a)} Gaussian waves, \textbf{(b)} triangle waves, \textbf{(c)} sawtooth waves and \textbf{(d)} pulse waves which are not found in the training dataset.}
	\label{Figure: 3}
\end{figure*}
We use the network architecture shown in Figure \ref{Figure: 1} to solve the Fisher's equation.The neural network architecture's input and output layers are synchronized with the spatial discretization, containing $n=128$ neurons each. The training utilized a global batch size of 32 and employed TensorFlow's mirrored strategy.
\par 
The first block $\phi$ is a transformer block containing multihead attention layer with $32$ heads, followed by layer normalization. Following layer normalization, there is a feedforward neural network consisting of two dense layers. These dense layers are equipped with rectified linear unit (ReLU) activation functions. After the feedforward neural network, the final layer in $\phi$ block consists of a layer normalization step.
\par
We trained the network on each of the types of initial conditions mentioned in the training dataset. This ensured that the network learned to generalize across different initial conditions and capture their underlying dynamics. After completing the training phase, we extensively tested the trained network using the test data, comprising trajectories with seven distinct types of initial conditions. Notably, the last four types of initial conditions: Gaussian waves, triangle waves, sawtooth waves, and pulse waves, were not encountered during the training process. In each scenario, the network prediction closely aligns with the exact solution, demonstrating excellent agreement between the two.
\par 
The inner encoder $\eta$, serves to diagonalize the system, and it also accomplishes the task of reducing the dimensionality of the system. The layer $v_k$ features a width distinct from that of the input and output layers. This disparity in width enables the dimensionality reduction within the $v_k$ layer, resulting in a low-dimensional model of rank $r$, we have opted for $r=21$. We utilize a single fully-connected linear layer for both the inner encoder $\eta$ and the inner decoder $\eta^{-1}$.
\par
\begin{figure*}
	\centering
	\includegraphics[width=0.80\columnwidth]{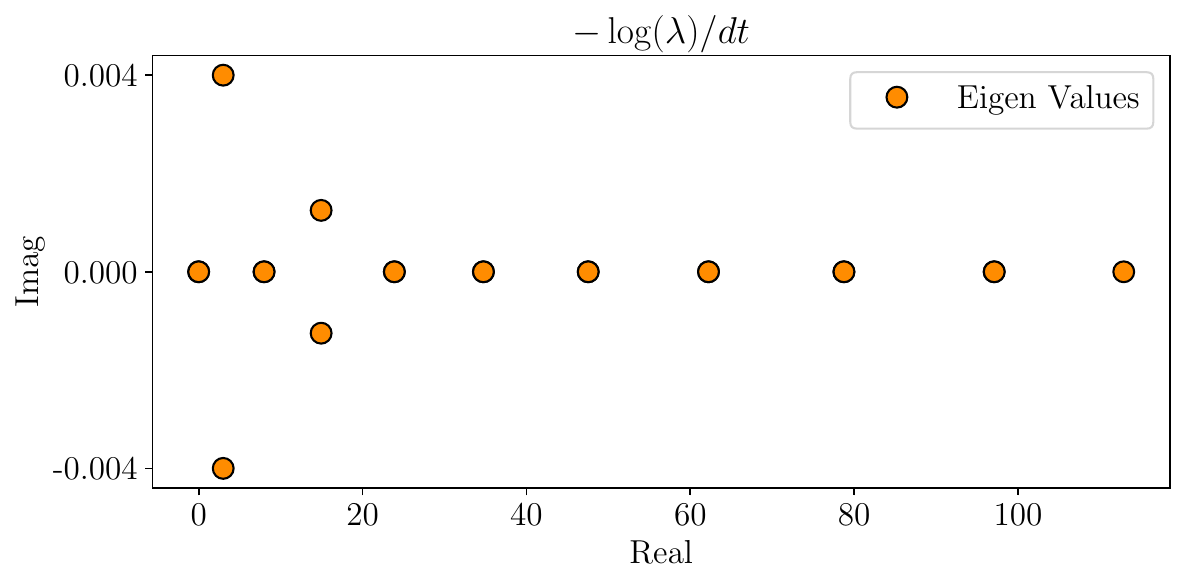}
	\caption{Eigen values of the approximated Koopman operator, $\textit{\textbf{K}}$.}
	\label{Figure: 4}
\end{figure*}
The comparison between the exact solution and the neural network predictions for Fisher's equation is presented in figure \ref{Figure: 2} and \ref{Figure: 3}. The first column in figure \ref{Figure: 2} and \ref{Figure: 3} plots the exact solution obtained by solving equation \ref{eq: 5} using finite difference technique. The second column in figure \ref{Figure: 2} and \ref{Figure: 3} shows the solution obtained using the proposed architecture. In figure \ref{Figure: 2}, the initial conditions represented are white noise, sine waves, and square waves. In figure \ref{Figure: 3}, the initial conditions represented are Gaussian waves, triangle waves, sawtooth waves, and pulse waves. It shall be noted that, the red line represents the initial condition, while the blue line shows the state of the wave at the final time, $t = 0.10$ seconds. In figure \ref{Figure: 2}, the initial conditions employed were of a similar type as provided to the architecture during the training phase. In figure \ref{Figure: 3}, the initial conditions utilized were not presented to the architecture during the training phase.  The comparison demonstrates that even with unseen data, our model exhibits strong performance, indicating its generalization across various types of initial conditions. The eigen values of the approximated Koopman operator, \textit{\textbf{K}}, has been shown in figure \ref{Figure: 4}. \par
\begin{figure*}
	\centering
	\includegraphics[width=0.70\columnwidth]{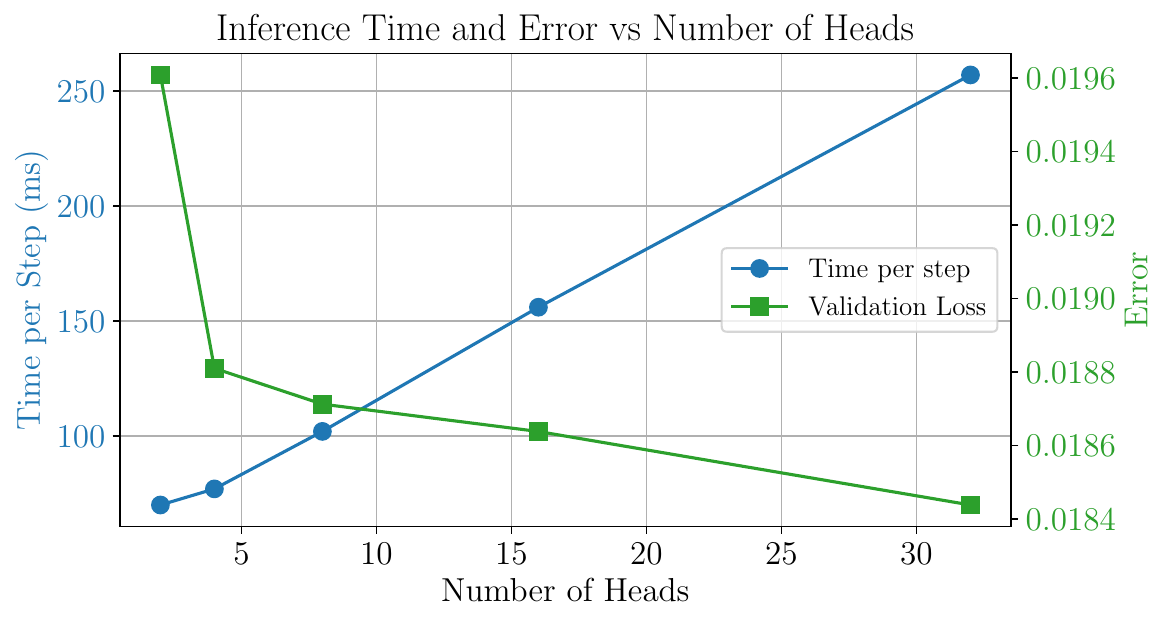}
	\caption{Relationship between inference time and error as a function of the number of attention heads.}
	\label{Figure: 5}
\end{figure*}

 Figure \ref{Figure: 5} illustrates the relationship between the number of attention heads and two critical metrics: inference time and error. The primary y-axis (in blue) shows that the inference time increases with the number of heads, ranging from approximately 70 ms at 2 heads to 257 ms at 32 heads. This trend suggests that the higher the number of attention heads, the greater the computational complexity and longer the processing times. The secondary y-axis (in green) shows that the total loss of the model decreases as the number of heads increases. This indicates that while adding more heads improves model accuracy, it comes at the cost of longer inference times and increased computational cost. The trade-off between these metrics depends on the specific application and user requirements, making it important to find an optimal balance between accuracy and computational efficiency. Furthermore, as the number of heads increases, the model’s scalability becomes a concern, as larger models require more resources and longer processing times. This can make models impractical for real-time or resource-limited environments.

\subsection{Comparative Analysis}
We conducted a comprehensive performance evaluation of our transformer-based outer encoder/decoder architecture against two alternative encoder/decoder structures as described by Gin et al. \cite{gin2021deep}: a dense encoder/decoder (DenseRes Block) and a convolutional-based encoder/decoder (ConvRes Block). The DenseRes Block comprised four densely connected layers with ReLU activation functions. The ConvRes Block comprises a series of convolutional layers followed by average pooling layers. It starts with 8 filters, progressively increasing to 64. All layers use a kernel size of 4, stride length of 1, and ReLU activation. Pooling layers have a size of 2 with no padding.
\par
\begin{figure}[t!]
	\begin{subfigure}{0.49\columnwidth}
		\centering
		\includegraphics[width=0.9\columnwidth]{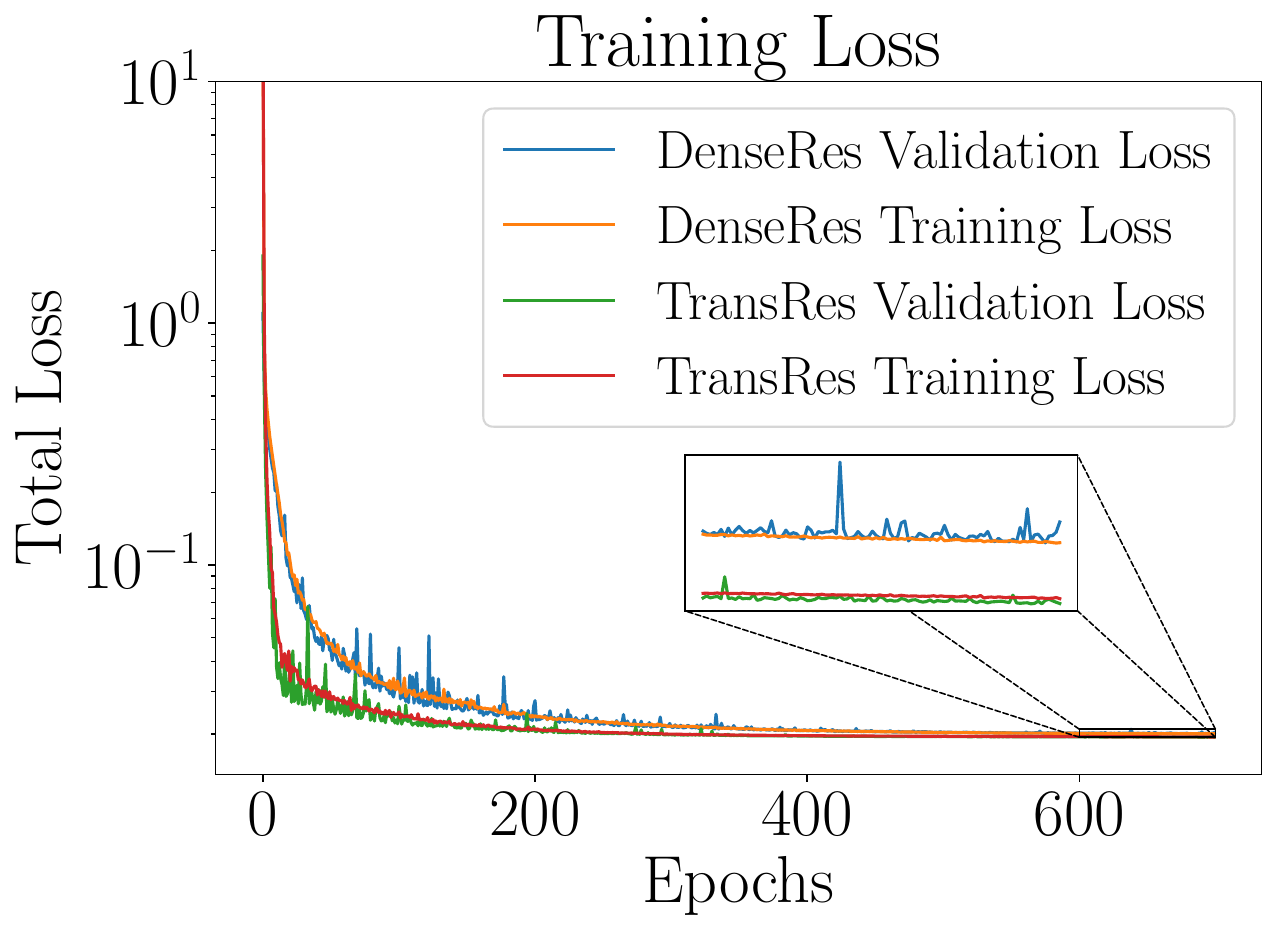}%
		\caption{}
		\label{Figure: 6_a}
	\end{subfigure}
	\begin{subfigure}{0.49\columnwidth}
		\centering
		\includegraphics[width=0.9\columnwidth]{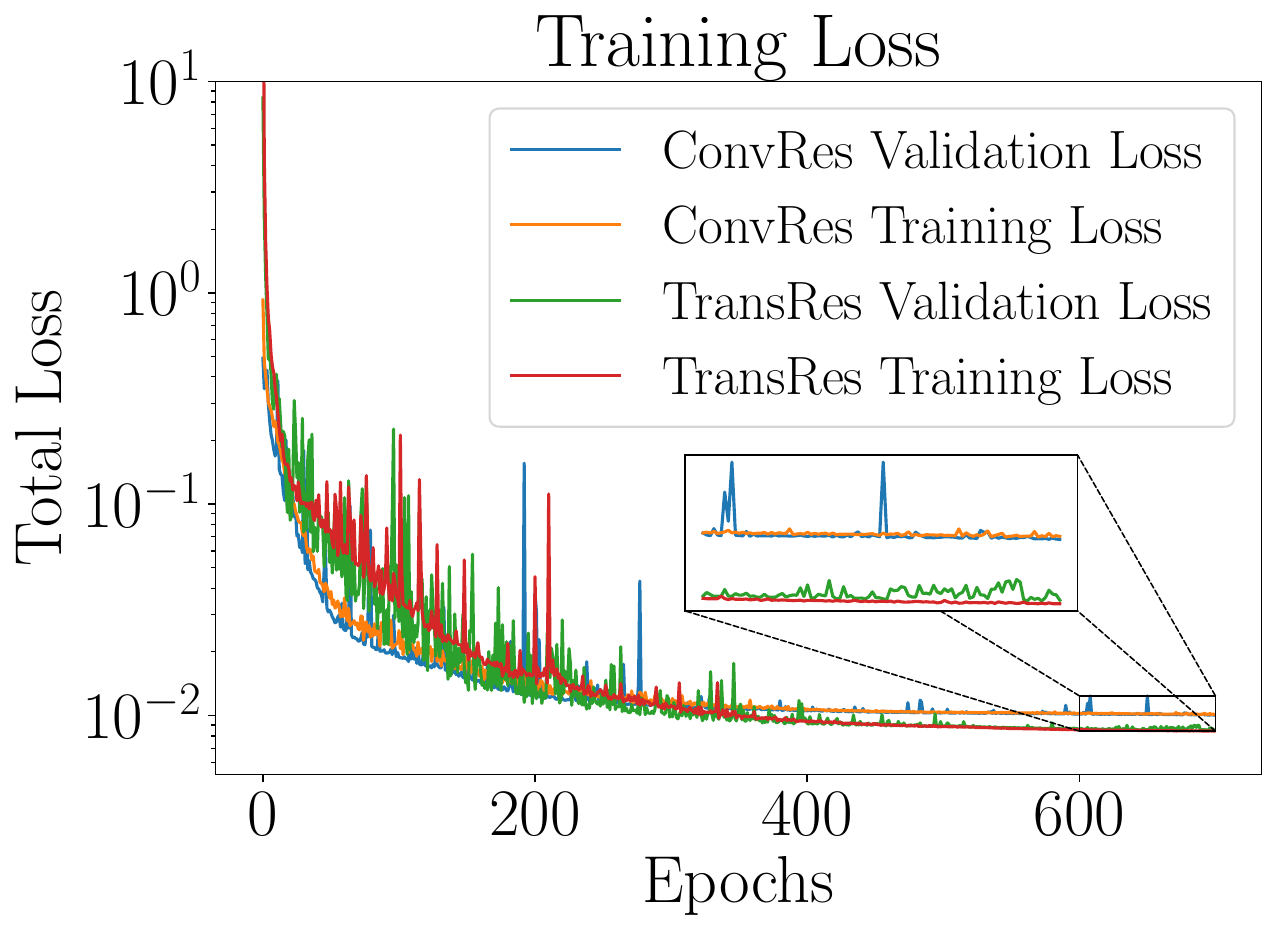}%
		\caption{}
		\label{Figure: 6_b}
	\end{subfigure}
	\caption{\textbf{(a)} Loss comparison of the TransRes Block with the DenseRes Block. \textbf{(b)} Loss comparison of the TransRes Block with the ConvRes Block. Both comparisons are presented using semilog plots. The inset plots illustrate the loss curves for the last 100 epochs.}
	\label{Figure: 6}
\end{figure}

The comparison with DenseRes Block was conducted using the Burger's equation dataset, comprising 60,000 initial conditions. It's evident from figure \ref{Figure: 6_a}, that both TransRes Block and DenseRes Block converge, with TransRes Block exhibiting slightly lower training and validation losses. Additionally, the loss curves indicate that TransRes Block achieves convergence in fewer epochs compared to DenseRes Block.
\par 
In the comparison with ConvRes Block, conducted using the Kuramoto-Shivashinisky equation dataset consisting of 60,000 initial conditions, a notable improvement in both training and validation losses is observed with TransRes Block, as shown in figure \ref{Figure: 6_b}.  While the loss curve suggests that ConvRes Block exhibited faster convergence during the initial epochs, it's noteworthy that TransRes Block eventually outperformed ConvRes Block,  resulting in superior performance and results.
\par
Our network gives better results in terms of the loss function for both Burger's and Kuramoto-Sivashinsky equation compared to Gin et al.'s \cite{gin2021deep}. They used two separate encoder blocks to solve two different systems. Despite using two separate encoder blocks, their results were not as par with the results, that we have obtained in this work. We are able to achieve better results even for three different types of systems using just one type of architecture. 

\par 
\begin{table*}[h]
	\caption{Loss metrics for the TransRes Block and DenseRes Block were evaluated on Burger's equation: $u_{t} + 10uu_{x} = u_{xx}$ over the domain $x \in (-\pi,\;\pi)$, subject to periodic boundary conditions.}
	\label{table: 1}
	\begin{tabularx}{\columnwidth}{
			|>{\centering\arraybackslash}X 
			|>{\centering\arraybackslash}X 
			|>{\centering\arraybackslash}X | }
		\hline
		\textbf{Outer Encoder/Decoder}$(\phi)$  & \textbf{Training Loss}  & \textbf{Validation Loss}\\
		\hline
		TransRes Block  & $1.944\times10^{-2}$  & $1.939\times10^{-2}$ \\
		DenseRes Block  & $2.003\times10^{-2}$& $2.036\times10^{-2}$ \\ 
		\hline
	\end{tabularx}
\end{table*}

\begin{table*}[h]
	\caption{Loss metrics for the TransRes Block and ConvRes Block were evaluated on Kuramoto-Sivashinsky equation: $u_{t} + u_{xx} + u_{xxxx} + uu_{x}= 0$ over the domain $x \in (-4\pi,\;4\pi)$, subject to periodic boundary conditions.}
	\label{table: 2}
	\begin{tabularx}{\columnwidth}{
			|>{\centering\arraybackslash}X 
			|>{\centering\arraybackslash}X 
			|>{\centering\arraybackslash}X | }
		\hline
		\textbf{Outer Encoder/Decoder}$(\phi)$ & \textbf{Training Loss}  & \textbf{Validation Loss}\\
		\hline
		TransRes Block  &  $8.433\times10^{-3}$& $8.516\times10^{-3}$ \\
		ConvRes Block  &  $1.012\times10^{-2}$&   $1.004\times10^{-2}$ \\
		\hline
	\end{tabularx}
\end{table*}
In table \ref{table: 1} and \ref{table: 2}, we have presented the values of total loss of TransRes Block, DenseRes Block and ConvRes Block, at 700$^{th}$ epoch. We can clearly see from table \ref{table: 1} that the value of training loss of our proposed network is relatively smaller than that of DenseRes Block. From table \ref{table: 2}, we can see that both the training and validation losses of our network are significantly lower than that of ConvRes's Block. 

\section{Discussion \& Conclusion}\label{Section: 4}
This paper introduces a deep learning architecture meticulously crafted to linearize Fisher's reaction-diffusion equation. Using Koopman operator theory, we successfully found coordinate transformations for the PDE, resulting in a low-dimensional model. This architecture's performance was then evaluated on additional PDEs, revealing its capacity to solve complex equations other than Fisher's equation with remarkable success. Our approach's consistent performance across different PDEs suggests that it could be used to solve a wide range of nonlinear problems. The architecture is scalable, and by increasing the number of heads or stacking more transformer layers, it can effectively linearize even highly nonlinear PDEs. It should be noted that our strategy is entirely data-driven, which eliminates the need for knowledge of the underlying equations. This makes it ideal for data where the governing equations are unknown.
\par 
In the future, further advancements can be made by optimizing these architectures. This could include not only changing the architecture and adjusting the loss functions, but also incorporating new advances in deep learning. Furthermore, continued research efforts might focus on improving the interpretability of the results produced by these frameworks. The application of these architectures extends beyond theoretical considerations to practical real-life problems. PDEs govern a variety of areas, including prey-predator interactions, disease modeling, climate dynamics, and many others, making these systems useful tools for understanding real-world occurrences.

\bibliographystyle{unsrt}  
\bibliography{references}  

%
%
%
%

\end{document}